\documentclass[]{article}

\usepackage[margin=1.5in]{geometry}  % set the margins to 1in on all sides
\usepackage{graphicx}              % to include figures
\usepackage{amsmath}               % great math stuff
\usepackage{amsfonts}              %for blackboard bold, etc
\usepackage{amsthm}                % better theorem environments
\usepackage{fancyhdr}
\usepackage{color}
\usepackage{hyperref}

\newtheorem{thm}{Theorem}[section]

\theoremstyle{definition}
\newtheorem{definition}{Definition}

\DeclareMathOperator{\pr}{pr}
\DeclareMathOperator{\Spec}{\mathrm{Spec}}

\title{Gaiotto's Lagrangian Subvarieties via Loop Groups} % \thanks{{\bf Preliminary version, not for distribution.}}
\author{Li, Yu}

\begin{document}

\maketitle

\begin{abstract}
The purpose of this note is to give a simple proof of the fact that a certain substack, defined in \cite{Ga}, of the moduli stack $T^{\ast}Bun_G(\Sigma)$ of Higgs bundles over a curve $\Sigma$, for a connected, simply connected semisimple group $G$, possesses a Lagrangian structure.  The substack, roughly speaking, consists of images under the moment map of global sections of principal $G$-bundles over $\Sigma$ twisted by a smooth symplectic variety with a Hamiltonian $G$-action.
\end{abstract}

\section{Introduction}

In this note we work over the complex numbers $\mathbb C$.  Throughout, let $G$ be a fixed connected, simply connected semisimple algebraic group and $\Sigma$ a fixed smooth irreducible projective curve, with a choice $K^{1/2}$ of square root of its canonical bundle.

In the differential-geometric setting, Gaiotto constructed in \cite{Ga}, for each symplectic vector space with a linear symplectic $G$-action, a subspace of the moduli space of stable Higgs bundles over $\Sigma$ and proved that this subspace is Lagrangian (see Appendix A of \cite{Ga}).  A more detailed proof in this setting, as well as explanations of the word `Lagrangian', can be found in Section 2 of \cite{H}.  In the algebro-geometric setting, a proof has recently been given in \cite{GR} (Theorem 1.3), based on the powerful theory of derived symplectic geometry.  The purpose of this note is to give an algebro-geometric proof of Gaiotto's statement in a down-to-earth manner.  In particular, our proof does not make use of any notions from derived algebraic geometry.  One interesting feature of our approach is that it allows us to see the picture very concretely - Gaiotto's statement boils down to the basic fact that the sum of residues of a rational $1$-form on $\Sigma$ equals zero.

Before stating the main theorem, let us fix some notations.  Let $Bun_G(\Sigma)$, or simply $Bun_G$, be the moduli stack of $G$-torsors over $\Sigma$.  It is well-known that the cotangent stack $T^{\ast}Bun_G$ is the moduli stack of Higgs bundles over $\Sigma$.  The stack $T^{\ast}Bun_G$, as a cotangent stack, has a canonical symplectic form $\Omega$ (see Definition \ref{def 2-form} below), defined as the differential of the Liouville $1$-form $\Lambda$ (see Definition \ref{def 1-form} below) on $T^{\ast}Bun_G$.  In the stacky context, there have been attempts to generalize the usual notions of isotropic/Lagrangian submanifolds to the notions of isotropic/Lagrangian substacks.  The interested readers are referred to page 1 of \cite{G} for more details.  In the derived context, such generalizations have been proposed in \cite{PTVV} (Definitions 2.7, 2.8 in \cite{PTVV}).  See also the definitions in Sections 2 and 3, as well as Section 5.2, of \cite{GR} for more discussions about these notions.  The main theorem of this note will be to the effect that a certain substack of $T^{\ast}Bun_G$ is Lagrangian. %(a certain morphism to $T^{\ast}Bun_G$ possesses a Lagrangian structure).

{\bf Remark.}  If it helps psychologically, one may replace $T^{\ast}Bun_G$ with $T^{\ast}Bun_G^{reg}$, the open substack of Higgs bundles without non-scalar automorphisms, so that all the usual notions of symplectic geometry make perfect sense.  Essentially nothing is lost by doing this.

Let $(X,\omega)$
%\footnote [1] {Here, we consider the more general case of symplectic schemes, as opposed to the case of a symplectic vector spaces as in \cite{Ga}.}
be a smooth symplectic variety with a Hamiltonian $G$-action preserving the symplectic form $\omega$ whose moment map is $\mu: X \rightarrow \mathfrak g^{\ast}$.  Assume moreover that there is a $\mathbb G_m$-action on $X$ such that the symplectic form $\omega$ has weight $2$ and the $\mathbb G_m$-action commutes with the $G$-action.  Given a $G$-torsor $P$, we write $X_P$ for $P \times _G X$.  If $(K^{1/2})^{\times}$ stands for the total space of the bundle $K^{1/2}$ without the zero section, then we write, by slightly abusing notation, $X_P \otimes K^{1/2}$ for $(X_P \times _{\Sigma} (K^{1/2})^{\times}) / \mathbb G_m$, where $\mathbb G_m$ acts on $X_P$ via its action on $X$, on $(K^{1/2})^{\times}$ by dilation, and on $X_P \times _{\Sigma} (K^{1/2})^{\times}$ by the diagonal action.

For the ease of future references, let us write $Y$ for the stack that parametrizes pairs $(P,s)$, where $P$ is a $G$-torsor over $\Sigma$ and $s$ is a global section of $X_P \otimes K^{1/2}$.  More precisely, given a $\mathbb C$-scheme $S$, the set of $S$-points of the stack $Y$ is defined as $$Y(S) := \{ (P,s): P \, \text{is a} \, G\text{-torsor over} \, \Sigma \times S, s \, \text{is a global section of} \, X_P \otimes \pr_1^{\ast} K^{1/2} \},$$ where $\pr_1: \Sigma \times S \rightarrow \Sigma$ is the natural projection to the first factor and $X_P$, as well as $X_P \otimes \pr_1^{\ast} K^{1/2}$, are constructed in an analogous way as in the previous paragraph.

Note that the moment map $\mu: X \rightarrow \mathfrak g^{\ast}$, being a homogeneous polynomial of degree $2$, induces, for each $\mathbb C$-scheme $S$ and each $G$-torsor $P$ over $\Sigma \times S$, a morphism
\begin{align} \label{mu}
X_P \otimes \pr_1^{\ast} K^{1/2} \rightarrow \mathfrak g^{\ast}_P \otimes \pr_1^{\ast} K,
\end{align}
where $\mathfrak g^{\ast}_P$ is the coadjoint bundle associated to $P$ and $K$ is the canonical bundle of $\Sigma$.
%Motivated by this, we define a stack $Z$ whose set of $S$-valued points, for a $\mathbb C$-scheme $S$, is
%\begin{align*}
%Z(S) := \{ (P, \phi): & P \, \text{is a} \, G\text{-torsor over} \, \Sigma \times S, \\ & \phi \in \Gamma(\Sigma \times S, \mathfrak g^{\ast}_P \otimes \pr_1^{\ast} K) \, \text{factors through the morphism} \, (\ref{mu}) \, \text{above}\}.
%\end{align*}
%It is then easy to see that the stack $Z$ is a substack of $T^{\ast}Bun_G$,
Since $T^{\ast}Bun_G$ parametrizes pairs $(P, \phi)$, where $P$ is a $G$-torsor over $\Sigma$ and $\phi$ is a global section of $\mathfrak g^{\ast}_P \otimes K$, the morphism (\ref{mu}) above induces a morphism $\tilde{\mu}: Y \rightarrow T^{\ast}Bun_G$.
%By abuse of notation, we also write $\tilde{\mu}$ for the composition $Y \rightarrow Z \hookrightarrow T^{\ast}Bun_G$.
The main theorem of this note is

%Since the moment map $\mu$ is a homogeneous polynomial of degree $2$, a section of the bundle $X_P \otimes K^{1/2}$ is mapped by $\mu$ to a section of the bundle $\mathfrak g^{\ast}_P \otimes K$, where $\mathfrak g^{\ast}_P$ is the coadjoint bundle associated to $P$ and $K$ is the canonical bundle of $\Sigma$.  Notice that $T^{\ast}Bun_G$ is the moduli stack parameterizing pairs $(P, \phi)$, where $\phi$ is a section of $\mathfrak g^{\ast}_P \otimes K$, by collecting the image of the morphism described above, as $P$ and the section $s$ of $X_P \otimes K^{1/2}$ vary, we obtain a substack $Z$ of $T^{\ast}Bun_G$.  Our notation for the morphism between $Y$ and $Z$ (or $T^{\ast}Bun_G$) will be $\tilde{\mu}$.

\begin{thm} \label{thm 1}
%The morphism $\tilde{\mu}: Y \rightarrow T^{\ast}Bun_G$ possesses a Lagrangian structure.  More precisely, t
The $2$-form $\tilde{\mu}^{\ast}(\Omega)$ vanishes on $Y$. %the preimage of $T^{\ast}Bun_G^{reg}$, the open substack of $T^{\ast}Bun_G$ corresponding to the Higgs bundles over $\Sigma$ without non-scalar automorphisms.
\end{thm}

{\bf Remark 1.}  If we write $Z$ for the substack of $T^{\ast}Bun_G$ whose set of $S$-valued points, for a $\mathbb C$-scheme $S$, is
\begin{align*}
Z(S) := \{ (P, \phi): & P \, \text{is a} \, G\text{-torsor over} \, \Sigma \times S, \\ & \phi \in \Gamma(\Sigma \times S, \mathfrak g^{\ast}_P \otimes \pr_1^{\ast} K) \, \text{factors through the morphism} \, (\ref{mu}) \, \text{above}\},
\end{align*}
then this theorem can be viewed as saying that the substack $Z$ of $T^{\ast}Bun_G$ is isotropic.

{\bf Remark 2.}  Our theorem is more flexible than Gaiotto's in \cite{Ga} in that we allow $X$ to be a symplectic variety with a Hamiltonian $G$-action, not just a symplectic vector space equipped with a linear symplectic $G$-action.

%{\bf Remark 3.}  For the sense in which the Theorem above can be interpreted as a Lagrangian condition, see the Remark in Section $2$ of \cite{H}.  More precise meaning of the word `Lagrangian', in the context of derived symplectic geometry, can be found in Definition 2.8 of \cite{GR}.  See Theorem 1.3 of \cite{GR} for the derived analog of Theorem \ref{thm 1}.

\section{The Proof}

In this section, we will try to be as humble and concrete as possible - whenever it helps with geometric intuition, we will pretend that $Bun_G$ and $T^{\ast}Bun_G$ are actual varieties, rather than stacks, and regard, for example, $Bun_G$ as consisting of ($\mathbb C$-valued) points given by transition morphisms with respect to some (fixed) cover of $\Sigma$.  The reason for this choice of perspective is that it is not hard to formalize what we are about to say rigorously in the stacky language (for a genaral $\mathbb C$-scheme $S$ and a $G$-torsor $P$ over $S$, apply Theorem 3 of \cite{DS} to trivialize $P$, then one can define the forms $\Lambda$ and $\Omega$ in an analogous manner to Definitions \ref{def 1-form} and \ref{def 2-form} and carry out essentially the same computation as in Sections \ref{pull-back} and \ref{genral-case} below), and that, probably more importantly, the readers will benefit from seeing the argument in concrete geometric terms.

Throughout, when we speak of tangent vectors to stacks, we mean {\it vecteur tangent universel} in the sense of \cite{LMB} (see page 169).  For the meaning of cotangent vectors we refer the readers to Th\'eor\`eme 17.16 of \cite{LMB}.

The plan of our proof is as follows.  We will first give a formula for the Liouville $1$-form $\Lambda$ on $T^{\ast}Bun_G$, then differentiate it to obtain a formula for the canonical symplectic form $\Omega$.  We then pull the form $\Omega$ back to the stack $Y$ and argue that the pull-back equals zero.

\subsection{Trivializations} \label{triv}

Let $P$ be a given $G$-torsor over $\Sigma$.  We first consider the special case where $P$, as well as $K^{1/2}$, can be trivialized over $\Sigma^{\circ} := \Sigma - \{\text{pt}\}$ for some point $\text{pt} \in \Sigma$ (see Satz 3.3 of \cite{Ha}), although this is in general not always the case.  The general case will follow easily as long as we fully understand this special case.

Let $D := \Spec \mathbb C[[z]]$ be a formal disk centered at $\text{pt}$.  Since any bundle over $D$ is trivial, $P$ and $K^{1/2}$ are determined by transition morphisms.
%Then the curve $\Sigma$ can be covered by the two charts $\Sigma^{\circ}$ and $D$.
Define $D^{\times} := \Spec \mathbb C((z))$ to be the punctured formal disk.  Then we write $g^{-1} \in  G(D^{\times})$ (resp. $T: D^{\times} \rightarrow \mathbb G_m$) for the transition morphism from $\Sigma^{\circ}$ to $D$ for $P$ (resp. $K^{1/2}$).  More precisely, for closed points, $(z,h) \in D^{\times} \times G$ over $\Sigma^{\circ}$ is identified with $(z, g(z)^{-1}h) \in D^{\times} \times G$ over $D$.  Similarly, $(z,t) \in D^{\times} \times \mathbb A^1$ over $\Sigma^{\circ}$ is identified with $(z, T(z)t) \in D^{\times} \times \mathbb A^1$ over $D$.

With these notations, the bundle $X_P \otimes K^{1/2}$ has a simple description as follows.  Over $\Sigma^{\circ}$, since $P$ is $G$-equivariantly isomorphic to the trivial $G$-torsor, $X_P$ is isomorphic to $\Sigma^{\circ} \times X$.  Since $K^{1/2}$ is $\mathbb G_m$-equivariantly isomorphic to $\Sigma^{\circ} \times \mathbb G_m$ over the same chart, we have the isomorphism $X_P \otimes K^{1/2} \simeq \Sigma^{\circ} \times X$ over $\Sigma^{\circ}$.  In a similar way, we see that $X_P \otimes K^{1/2} \simeq D \times X$ over $D$.  Using the transition morphisms from the previous paragraph, we see that the transition morphism for $X_P \otimes K^{1/2}$ from $\Sigma^{\circ}$ to $D$ is $T^{-1}g^{-1}$.

\subsection{The `Source Stack' $Y$} \label{sstk}

Let us describe the source stack $Y$ in this section.
%\footnote {See our discussion prior to \ref{thm 1}.}
In concrete words, $Y$ parametrizes equivalence classes of triples $(g, s^{\circ}, s') \in G(D^{\times}) \times Map(\Sigma^{\circ}, X) \times Map(D,X)$ satisfying the equation 
\begin{align} \label{prel}
T^{-1}g^{-1}s^{\circ} = s',
\end{align}
where $(g, s^{\circ}, s')$ is equivalent to $(h, t^{\circ}, t')$ if there exist $g_1 \in G(\Sigma^{\circ})$ and $g_2 \in G(D)$ such that $g_1gg_2^{-1} = h$, $g_1s^{\circ} = t^{\circ}$, and $g_2^{-1}s' = t'$.  In the sequel, for simplicity, we will suppress the words `equivalence classes of' and write only `triples' or `tuples' in all similar situations.
%We sometimes also think of $g$ as a morphism from $D^{\times}$ to $G$, although it is in fact a double coset of the loop group $G(D^{\times})$.

Given a point $(g, s^{\circ}, s') \in Y$, $T_{(g, s^{\circ}, s')}Y$ consists of $1$st order infinitesimal deformations of $(g, s^{\circ}, s')$.  Concretely, a tangent vector to $Y$ at $(g, s^{\circ}, s')$ consists of a triple $(\dot g, \dot s^{\circ}, \dot {s}') \in \mathfrak g(D^{\times}) \times Map(\Sigma^{\circ}, TX) \times Map(D, TX)$ so that $\pi \circ \dot s^{\circ} = s^{\circ}$, $\pi \circ \dot{s}' = s'$, where $\pi: TX \rightarrow X$ is the projection, and
\begin{align} \label{rel}
\dot{s}' = T^{-1}g^{-1}\dot s^{\circ} - \rho_{s'}(\dot g),
\end{align}
where $\rho_{s'}(\dot g)$ stands for the infinitesimal action of $\dot g$ at $s'$.  Note that equation (\ref{rel}) is obtained from equation (\ref{prel}) by considering $1$st order infinitesimal deformations.

Heuristically, one can view tangent vectors as infinitesimal curves and write $s^{\circ}_t$ (resp. $s'_t$) for the infinitesimal curve representing $\dot s^{\circ}$ (resp. $\dot{s}'$).  Using the suggestive notation $g\exp(t \dot g)$ for an infinitesimal curve representing $\dot g$, we see from (\ref{prel}) above that we must have the equation $T^{-1}(g\exp(t\dot g))^{-1} s^{\circ}_t = s'_t$.  Applying $\frac {\,d}{\,dt}|_{t=0}$ and using the Leibniz rule, we arrive at equation (\ref{rel}).

{\bf Remark 1.}  This heuristic argument can be made precise by using the dual numbers.

{\bf Remark 2.}  %In the heuristic argument we are thinking of $g$ as an element of the loop group $G(D^{\times})$ and, hence, $\dot g$ as an element of the loop Lie algebra $\mathfrak g(D^{\times})$.  We will make use of these two view points interchangeably, without making further remarks.
In the argument above we have tacitly used the assumption that all $G$-torsors $Q$ sufficiently close to $P$ in $Bun_G$ can be trivialized over $\Sigma^{\circ}$.  See Theorem 3 in \cite{DS} for more details about this.

\subsection{The Liouville $1$-form $\Lambda$ on $T^{\ast}Bun_G$} \label{1-form}

Recall that, given a $G$-torsor $P$, one has $T^{\ast}_PBun_G = H^0(\Sigma, \mathfrak g^{\ast}_P \otimes K)$, where $\mathfrak g^{\ast}_P$ stands for the coadjoint bundle associated to $P$ and $K$ stands for the canonical bundle of $\Sigma$.  Since $T$ is the transition morphism for $K^{1/2}$, the transition morphism for $K$ is given by $T^2$.  The same analysis as in Section \ref{triv} then tells us that a global section of $\mathfrak g^{\ast}_P \otimes K$ consists of a pair $(\phi^{\circ}, \phi') \in \mathfrak g^{\ast}(\Sigma^{\circ}) \times \mathfrak g^{\ast}(D)$ satisfying $\phi' = T^{-2}_{\cdot} g_{\cdot} \phi^{\circ}$, where $g \in G(D^{\times})$ corresponds to $P$ and acts on $\phi^{\circ}$ by the coadjoint action.

Given a triple $(g, \phi^{\circ}, \phi')$ as in the previous paragraph, by the same argument as in Section \ref{sstk} we see that a tangent vector to $T^{\ast}Bun_G$ at $(g, \phi^{\circ}, \phi')$ consists of a triple $(\dot g, \dot{\phi^{\circ}}, \dot{\phi}') \in \mathfrak g(D^{\times}) \times \mathfrak g^{\ast}(\Sigma^{\circ}) \times \mathfrak g^{\ast}(D)$ satisfying an equation similar to (\ref{rel}) in Section \ref{sstk}.  The push-forward to $Bun_G$ of this tangent vector is just $\dot g: D^{\times} \rightarrow \mathfrak g$.  With these in mind, it makes sense to make the following

\begin{definition} \label{def 1-form}
The Liouville $1$-form $\Lambda$ is given by $$\Lambda(\dot g, \dot \phi^{\circ}, \dot{\phi}') := (g, \phi^{\circ}, \phi')(\dot g) = Res_{z=0}(\langle \phi', \dot g \rangle \,dz),$$ where we used $\,dz$ to trivialize $K$ over $D$.
\end{definition}

\subsection{The Symplectic $2$-form $\Omega$ on $T^{\ast}Bun_G$} \label{2-form}

Given $(g, \phi^{\circ}, \phi') \in T^{\ast}Bun_G$ and tangent vectors $(\dot g_1, \dot \phi^{\circ}_1, \dot{\phi}'_1)$, $(\dot g_2, \dot \phi^{\circ}_2, \dot{\phi}'_2)$ to $T^{\ast}Bun_G$ at $(g, \phi^{\circ}, \phi')$, we define

\begin{definition} \label{def 2-form}
The symplectic $2$-form $\Omega$ on $T^{\ast}Bun_G$ is given by
\begin{align}
&\Omega((\dot g_1, \dot {\phi}^{\circ}_1, \dot{\phi}'_1), (\dot g_2, \dot {\phi}^{\circ}_2, \dot{\phi}'_2)) \nonumber \\
:=& Res_{z=0} (\langle \phi'_1, \dot g_2 \rangle \,dz) - Res_{z=0} (\langle \phi'_2, \dot g_1 \rangle \,dz) - Res_{z=0} (\langle \phi', [\dot g_1, \dot g_2] \rangle \,dz).
\end{align}
\end{definition}

This section is devoted to the justification of this definition.  We will interpret $\Omega$ as obtained from $\Lambda$ by Cartan's formula for exterior derivatives.

%In this section we compute $\Omega((\dot g_1, \dot \phi^{\circ}_1, \dot{\phi}'_1), (\dot g_2, \dot \phi^{\circ}_2, \dot{\phi}'_2))$ for tangent vectors $(\dot g_1, \dot \phi^{\circ}_1, \dot{\phi}'_1)$ and $(\dot g_2, \dot \phi^{\circ}_2, \dot{\phi}'_2)$ to $T^{\ast}Bun_G$ at $(g, \phi^{\circ}, \phi')$.  We will do this using Cartan's formula for exterior derivatives.
%the formula that, for a manifold $M$ and a $d$-form $\alpha$ on $M$, 

%\begin{align*}
%& d\alpha(\xi_0, \cdots, \xi_d) \\ = &\sum_{i=0}^d (-1)^i \xi_i \alpha(\xi_0, \cdots, \hat{\xi_i}, \cdots, \xi_d) + \sum_{i<j} (-1)^{i+j} \alpha ([\xi_i,\xi_j], \xi_0, \cdots, \hat{\xi_i}, \cdots, \hat{\xi_j}, \cdots, \xi_d),
%\end{align*}
%where $\xi_0, \cdots, \xi_d$ are vector fields on $M$.

Notice first that the $1$-form $\Lambda$ is actually defined on the space $\mathfrak g(D^{\times}) \times \mathfrak g^{\ast} (\Sigma^{\circ}) \times \mathfrak g^{\ast}(D)$,
%of all triples $(\dot g, \dot {\phi}^{\circ}, \dot{\phi}') \in TBun_G \times \mathfrak g^{\ast}(\Sigma) \times \mathfrak g^{\ast}(D)$,
not just the subspace $T^{\ast}Bun_G$ of $\mathfrak g(D^{\times}) \times \mathfrak g^{\ast} (\Sigma^{\circ}) \times \mathfrak g^{\ast}(D)$.  So, if we write $\tilde{\Lambda}$ for this $1$-form on the ambient space $\mathfrak g(D^{\times}) \times \mathfrak g^{\ast} (\Sigma^{\circ}) \times \mathfrak g^{\ast}(D)$ and $\iota$ for the inclusion of $T^{\ast}Bun_G$ into the ambient space, we have $$\,d \Lambda = \,d \iota^{\ast} \tilde{\Lambda} = \iota^{\ast} \,d \tilde{\Lambda}.$$

For $i=1, 2$, extend the tangent vectors $(\dot g_i, \dot \phi^{\circ}_i, \dot{\phi}'_i)$ to vector fields $(\dot G_i, \dot \Phi^{\circ}_i, \dot{\Phi}'_i)$ on $\mathfrak g(D^{\times}) \times \mathfrak g^{\ast} (\Sigma^{\circ}) \times \mathfrak g^{\ast}(D)$ near $(g, \phi^{\circ}, \phi')$.  Since in the computation of $\, d \tilde{\Lambda} ((\dot g_1, \dot \phi^{\circ}_1, \dot{\phi}'_1), (\dot g_2, \dot \phi^{\circ}_2, \dot{\phi}'_2))$, it does not matter which extension we use, we will choose the one where $\dot G_i(h,\psi^{\circ}, \psi') = \dot g_i$ for $i=1, 2$, where $(h, \psi^{\circ}, \psi')$ is a point in $\mathfrak g(D^{\times}) \times \mathfrak g^{\ast} (\Sigma^{\circ}) \times \mathfrak g^{\ast}(D)$ near $(g, \phi^{\circ}, \phi')$.

%The alerted readers may have doubts about the existence of such extensions, since we want the triples $(\dot G_i, \dot \Phi^{\circ}_i, \dot{\Phi}'_i)$ ($i=1, 2$) to be vector fields on $T^{\ast}Bun_G$, which is a subspace of $\mathfrak g(D^{\times}) \times \mathfrak g^{\ast} (\Sigma^{\circ}) \times \mathfrak g^{\ast}(D)$, rather than vector fields on the ambient space $\mathfrak g(D^{\times}) \times \mathfrak g^{\ast} (\Sigma^{\circ}) \times \mathfrak g^{\ast}(D)$ (see the last paragraph of Section \ref{1-form}).
%We argue that we do not need to worry about this issue.  In fact, the relation among the three components is the $1$st order infinitesimal deformation of the relation $\phi' = T^{-2}_{\cdot} g_{\cdot} \phi^{\circ}$ on $T^{\ast}Bun_G$.
%This computation tells us that, in order to compute $\, d \Lambda$, we can first pullback to the ambient space, where we do not need to worry about the existence of extensions, then compute exterior derivative, and then restrict back to $T^{\ast}Bun_G$.

%and that the pull-back along the inclusion of $T^{\ast}Bun_G$ into $F$ commutes with exterior differential.  As a consequence, we see that, in order to compute $\Omega$, it does not matter whether we compute it on $T^{\ast}Bun_G$ or on $F$.  Since there is no relations on the latter stack, there is no relations among the components of the tangent vectors to the latter stack.  Therefore, we do not need to worry about existence of extensions of tangent vectors to vector fields.

Observe that the function $\tilde{\Lambda}(\dot G_2, \dot \Phi^{\circ}_2, \dot {\Phi}'_2)$ sends a point $(h, \psi^{\circ}, \psi')$ near the point $(g, \phi^{\circ}, \phi')$ to $Res_{z=0}(\langle \psi', \dot g_2 \rangle \,dz)$.  Applying the vector field $(\dot G_1, \dot \Phi^{\circ}_1, \dot {\Phi}'_1)$ to the function $\tilde{\Lambda}(\dot G_2, \dot \Phi^{\circ}_2, \dot {\Phi}'_2)$ and evaluating at $(g, \phi^{\circ}, \phi')$, we get $$Res_{z=0}(\langle \dot{\phi}'_1, \dot g_2 \rangle \,dz).$$
%To see this more concretely, we let $(g\exp(t \dot g_1), \phi^{\circ}_{1,t}, \phi'_{1,t})$ be an infinitesimal curve representing $(\dot g_1, \dot \phi^{\circ}_1, \dot{\phi}'_1)$.  Then what we want is just $\frac {\,d}{\,dt}|_{t=0} Res_{z=0} (\langle \phi'_{1,t}, \dot g_2 \rangle \,dz)$, which is equal to $Res_{z=0} (\langle \dot {\phi}'_1, \dot g_2 \rangle \,dz)$.
Similarly, applying the vector field $(\dot G_2, \dot \Phi^{\circ}_2, \dot {\Phi}'_2)$ to the function $\tilde{\Lambda}(\dot G_1, \dot \Phi^{\circ}_1, \dot {\Phi}'_1)$ and evaluating at $(g, \phi^{\circ}, \phi')$, we get $$Res_{z=0}(\langle \dot{\phi}'_2, \dot g_1 \rangle \,dz).$$  Finally, observe that the push-forward of $[(\dot G_1, \dot {\Phi}^{\circ}_1, \dot{\Phi}'_1), (\dot G_2, \dot {\Phi}^{\circ}_2, \dot{\Phi}'_2)]$ to $Bun_G$ at $(g, \phi^{\circ}, \phi')$ is just $[\dot g_1, \dot g_2]$, so we see that
$$\tilde{\Lambda} ([(\dot G_1, \dot {\Phi}^{\circ}_1, \dot{\Phi}'_1), (\dot G_2, \dot {\Phi}^{\circ}_2, \dot{\Phi}'_2)])
= Res_{z=0} (\langle \phi', [\dot g_1, \dot g_2] \rangle \,dz)$$
at the point $(g, \phi^{\circ}, \phi')$.  Using Cartan's formula for exterior derivatives, we have
\begin{align} \label{eqn1}
&\, d \tilde{\Lambda}((\dot g_1, \dot {\phi}^{\circ}_1, \dot{\phi}'_1), (\dot g_2, \dot {\phi}^{\circ}_2, \dot{\phi}'_2)) \nonumber \\
=& Res_{z=0} (\langle \phi'_1, \dot g_2 \rangle \,dz) - Res_{z=0} (\langle \phi'_2, \dot g_1 \rangle \,dz) - Res_{z=0} (\langle \phi', [\dot g_1, \dot g_2] \rangle \,dz),
\end{align}
thus justifying Definition \ref{def 2-form}.

\subsection{Pull-back of $\Omega$ to $Y$} \label{pull-back}

We now put together everything we have seen so far.  Let $(g, s^{\circ}, s')$ be a $\mathbb C$-valued point in the source stack $Y$ and let $(\dot g_i, \dot s^{\circ}_i, \dot{s}'_i)$ be tangent vectors to $Y$ at $(g, s^{\circ}, s')$ ($i=1,2$).  These two vectors are pushed forward by $\tilde{\mu}$ to $(\dot g_i, \,d\mu (\dot s^{\circ}_i), \,d \mu(\dot{s}'_i))$ ($i=1,2$).  Using formula (\ref{eqn1}) we have seen in Section \ref{2-form}, we see that 

\begin{align*}
& \Omega((\dot g_1, \,d \mu(\dot s^{\circ}_1), \,d \mu (\dot{s}'_1)), (\dot g_2, \,d \mu(\dot s^{\circ}_2), \,d \mu (\dot{s}'_2))) \\
= & Res_{z=0} (\langle \, d\mu (\dot{s}'_1), \dot g_2 \rangle \,dz - \langle \, d\mu (\dot{s}'_2), \dot g_1 \rangle \,dz - \langle \mu \circ s', [\dot g_1, \dot g_2] \rangle \,dz) \\
= & Res_{z=0} (- \omega(\dot{s}'_1, \rho_{s'}(\dot g_2)) \,dz + \omega(\dot{s}'_2, \rho_{s'}(\dot g_1)) \,dz - \langle \mu \circ s', [\dot g_1, \dot g_2] \rangle \,dz).
\end{align*}

Now let $\alpha$ be the nowhere vanishing $1$-form on $\Sigma^{\circ}$ used in the trivialization of $K$.  So we have $\alpha = T^{-2}\,dz$ on $D^{\times}$.  Using equation (\ref{rel}) from Section \ref{sstk}, we have

\begin{align} \label{eqn2}
& \omega(g^{-1} \dot s^{\circ}_1, g^{-1} \dot s^{\circ}_2) \alpha \nonumber \\
= & \omega(g^{-1} \dot s^{\circ}_1, g^{-1} \dot s^{\circ}_2) T^{-2} \,dz \nonumber \\
= & \omega(T^{-1} g^{-1} \dot s^{\circ}_1, T^{-1} g^{-1} \dot s^{\circ}_2) \,dz \nonumber \\
= & \omega(\dot{s}'_1 + \rho_{s'}(\dot g_1), \dot{s}'_2 + \rho_{s'}(\dot g_2)) \,dz \nonumber \\
= & \omega(\dot{s}'_1, \dot{s}'_2) \,dz + \omega(\dot{s}'_1, \rho_{s'}(\dot g_2)) \,dz - \omega(\dot{s}'_2, \rho_{s'}(\dot g_1)) \,dz + \omega(\rho_{s'}(\dot g_1), \rho_{s'}(\dot g_2)) \,dz \nonumber \\
= & \omega(\dot{s}'_1, \dot{s}'_2) \,dz + \omega(\dot{s}'_1, \rho_{s'}(\dot g_2)) \,dz - \omega(\dot{s}'_2, \rho_{s'}(\dot g_1)) \,dz + \langle \mu \circ s', [\dot g_1, \dot g_2] \rangle \,dz
\end{align}

The very first line of (\ref{eqn2}), by $G$-invariance of $\omega$, is equal to $\omega(\dot{s}^{\circ}_1, \dot s^{\circ}_2) \alpha$, which is a $1$-form defined on $\Sigma$ with a pole at $\text{pt}$.  Since the sum of residues of a rational $1$-form on $\Sigma$ is zero, we see that the residue at $z=0$ of the bottom line of (\ref{eqn2}) is zero.  Since $\omega(\dot{s}'_1, \dot{s}'_2) \,dz$ is a $1$-form defined on the disk $D$, its residue at $z=0$ is also zero.  This implies that the residue at $z=0$ of the last three terms in the bottom line of (\ref{eqn2}) is zero, thus proving that $\tilde{\mu}^{\ast}\Omega$ vanishes. %on the preimage of $T^{\ast}Bun_G^{reg}$.

\subsection{The General Case} \label{genral-case}

Finally we remove the assumption that $P$ and $K^{1/2}$ can be trivialized over $\Sigma - \{\text{pt}\}$ as follows.  In this general case we remove finitely many points $\text{pt}_1, \cdots, \text{pt}_n$ so that $P$ and $K^{1/2}$ are can be trivialized over $\Sigma^{\circ} := \Sigma - \{ \text{pt}_1, \cdots, \text{pt}_n \}$ (see Satz 3.3 of \cite{Ha} and Theorem 3 of \cite{DS} for more details).  As before, we write $D_i := \Spec \mathbb C[[z]]$ ($i=1, \cdots ,n$) for formal disks centered at $\text{pt}_i$ and $D_i^{\times} := \Spec \mathbb C((z))$ for punctured formal disks centered at $\text{pt}_i$.  We write $g_i^{-1} \in G(D_i^{\times})$ ($i=1, \cdots, n$) for the transition morphisms of $P$ from $\Sigma^{\circ}$ to $D_i$ and $T_i: D_i^{\times} \rightarrow \mathbb G_m$ for the transition morphisms of $K^{1/2}$ from $\Sigma^{\circ}$ to $D_i$.  Using these transitions, the bundle $X_P \otimes K^{1/2}$ can be described by the transition morphisms $T_i^{-1}g_i^{-1}$ for $i=1, \cdots, n$, and, hence, the source stack $Y$ consists of tuples $(g_i, s^{\circ}, s'_i) \in \prod_{i=1}^n G(D^{\times}_i) \times Map(\Sigma^{\circ}, X) \times \prod_{i=1}^n Map(D_i,X)$ satisfying the equations $s'_i = T_i^{-1}g_i^{-1}s^{\circ}$ for all $i$.

Given a tuple $(g_i, s^{\circ}, s'_i)$ as in the previous paragraph, a tangent vector to $Y$ at $(g_i, s^{\circ}, s'_i)$ is a tuple $(\dot g_i, \dot s^{\circ}, \dot{s}'_i) \in \prod_{i=1}^n \mathfrak g(D^{\times}_i) \times Map(\Sigma^{\circ}, TX) \times \prod_{i=1}^n Map(D_i, TX)$ such that $\pi \circ \dot s^{\circ} = s^{\circ}$, $\pi \circ \dot{s}'_i = \hat s_i$ and $$\dot{s}'_i = T_i^{-1}g_i^{-1}\dot s^{\circ} - \rho_{s'_i}(\dot g_i)$$ for all $i$, where $\pi: TX \rightarrow X$ is the natural projection.

For a point $(g_i, \phi^{\circ}, \phi'_i) \in T^{\ast}Bun_G$, and a tangent vector $(\dot g_i, \dot{\phi}^{\circ}, \dot{\phi}'_i)$ to $T^{\ast}Bun_G$ at $(g_i, \phi^{\circ}, \phi'_i)$, we define $\Lambda(\dot g_i, \dot{\phi}^{\circ}, \dot{\phi}'_i) := \sum_{i=1}^n Res_{z=0} (\langle \phi'_i, \dot g_i \rangle \,dz)$.  The exact same computation as in Section \ref{2-form} justifies the following definition: $$\Omega ((\dot g_{i,1}, \dot{\phi}^{\circ}_1, \dot{\phi}_{i,1}^{'}), (\dot g_{i,2}, \dot{\phi}^{\circ}_2, \dot{\phi}_{i,2}^{'})) := \sum_{i=1}^n Res_{z=0} (\langle \dot{\phi}_{i,1}^{'}, \dot g_{i,2} \rangle \,dz - \langle \dot{\phi}_{i,2}^{'}, \dot g_{i,1} \rangle \,dz - \langle \phi'_i, [\dot g_{i,1}, \dot g_{i,2}] \rangle \,dz)$$ for tangent vectors $(\dot g_{i,j}, \dot{\phi}^{\circ}_j, \dot{\phi}_{i,j}^{'})$ ($j=1,2$) to $T^{\ast}Bun_G$ at $(g_i, \phi^{\circ}, \phi'_i)$.  This tells us that, if $(g_i, s^{\circ}, s'_i)$ is a point in $Y$ and $(\dot g_{i,j}, \dot s^{\circ}_j, \dot{s}'_{i,j})$ are tangential to $Y$ at $(g_i, s^{\circ}, s'_i)$ ($j=1,2$), then

\begin{align}  \label{eqn3}
& \Omega ((\dot g_{i,1}, \,d \mu (\dot s^{\circ}_1), \,d \mu (\dot{s}'_{i,1})), (\dot g_{i,2}, \,d \mu (\dot s^{\circ}_2), \,d \mu (\dot{s}'_{i,2}))) \nonumber \\
= & \sum_{i=1}^n Res_{z=0} (-\omega (\dot{s}'_{i,1}, \rho_{s'_i}(\dot g_{i,2})) \,dz + \omega (\dot{s}'_{i,2}, \rho_{s'_i}(\dot g_{i,1})) \,dz - \langle \mu \circ s'_i, [\dot g_{i,1}, \dot g_{i,2}] \rangle \,dz).
\end{align}

We conclude that the bottom line of (\ref{eqn3}) is zero.  Recall that $\alpha$ is a nowhere vanishing $1$-form on $\Sigma^{\circ}$ used in the trivialization of $K$.  A computation similar to the one in the previous section tells us that
\begin{align*}
&\omega (\dot{s}'_{i,1}, \dot{s}'_{i,2}) \, dz - \omega (\dot s^{\circ}_1, \dot s^{\circ}_2) \alpha \nonumber \\
=& -\omega (\dot{s}'_{i,1}, \rho_{s_i}'(\dot g_{i,2})) \,dz + \omega (\dot{s}'_{i,2}, \rho_{s'_i}(\dot g_{i,1})) \,dz - \langle \mu \circ s'_i, [\dot g_{1,i}, \dot g_{2,i}] \rangle \,dz
\end{align*}
for all $i$.  Since $\omega (\dot{s}'_{i,1}, \dot{s}'_{i,2}) \, dz$ is a $1$-form defined on $D_i$, its residue at $z=0$ is zero.  From this we see that the bottom line of (\ref{eqn3}) equals $-\sum_{i=1}^n Res_{z=0} (\omega (\dot s^{\circ}_1, \dot s^{\circ}_2)\alpha)$.  This is zero, again because the sum of residues of a rational $1$-form on $\Sigma$ is equal to zero.

{\bf Acknowledgments.}  The author would like to thank V. Ginzburg for introducing this problem to the author, for his continual support and for many invaluable discussions during the preparation of this note.  The author is also grateful to A. Beilinson for answering several questions about trivializations of principal $G$-bundles and to M. Nori for discussions about a preliminary version of this note.

\makebox
{\bf Li, Yu:} Department of Mathematics, University of Chicago, Chicago, IL 60637, USA; \par
{\bf liyu@math.uchicago.edu}

\end{document}